\documentclass[11pt]{amsart}
\usepackage{amsmath,amssymb,amsthm,amscd,verbatim}
\usepackage{graphicx}
\usepackage{epsfig}
\usepackage{hyperref}

\newtheorem{thm}{Theorem}
\newtheorem{lem}[thm]{Lemma}
\newtheorem{conj}[thm]{Conjecture}

\newtheorem{cor}[thm]{Corollary}
\newtheorem{cla}[thm]{Claim}
\newtheorem*{thm*}{Theorem}

\theoremstyle{definition}

\newtheorem*{lpc}{Lov\'asz-Plummer conjecture}

\theoremstyle{remark}

\def\ceps{3656}

\newcommand{\calY}{\mathcal{Y}}
\newcommand{\calM}{\mathcal{M}}

\newcommand{\calC}{\mathcal{C}}

\newcommand{\calX}{\mathcal{X}}

\newcommand{\Ex}{\mathbb{E}}

\newcommand{\Sone}{{\bf [S1] }}
\newcommand{\Stwo}{{\bf [S2] }}
\newcommand{\fw}{\mathit{fw}}
\newcommand{\bM}{\mathbf{M}}

\begin{document}
\title{Exponentially many perfect matchings in cubic graphs}
\author{Louis Esperet} \address{Laboratoire G-SCOP (CNRS,
  Grenoble-INP), Grenoble, France}
\email{louis.esperet@g-scop.fr}\thanks{The first author is partially
  supported by ANR Project Heredia under Contract
  \textsc{anr-10-jcjc-heredia}.}

\author{Franti\v{s}ek Kardo\v{s}} \address{Institute of Mathematics,
  Faculty of Science, Pavol Jozef \v{S}af\'arik University,
  Ko\v{s}ice, Slovakia} \email{frantisek.kardos@upjs.sk} \thanks{The
  second author is partially supported by Slovak Research and
  Development Agency under Contract No. APVV-0007-07.}

\author{Andrew D.\ King} \address{Department of Industrial Engineering
  and Operations Research, Columbia University, New-York, NY, USA}
\email{andrew.d.king@gmail.com} \thanks{The third author is supported
  by an NSERC Postdoctoral Fellowship.}

\author{Daniel Kr{\'a}l'} \address{Department of Applied Mathematics
  and Institute for Theoretical Computer Science, Faculty of
  Mathematics and Physics, Charles University, Malostransk{\'e}
  n{\'a}m{\v e}st{\'\i} 25, 118 00 Prague, Czech Republic.}
\email{kral@kam.mff.cuni.cz} \thanks{The Institute for Theoretical
  Computer Science is supported by Ministry of Education of Czech
  Republic as projects 1M0545. The work leading to this invention has
  received funding from the European Research Council under the
  European Union's Seventh Framework Programme (FP7/2007-2013)/ERC
  grant agreement no. 259385.}

\author{Serguei Norine} \address{Department of Mathematics, Princeton
  University, Princeton, NJ, USA}
\email{snorin@math.princeton.edu}\thanks{The fifth author is partially
  supported by NSF under Grant No. DMS-0803214.}

\date{}

\begin{abstract}
  We show that every cubic bridgeless graph $G$ has at least
  $2^{|V(G)|/\ceps}$ perfect matchings. This confirms an old
  conjecture of Lov\'{a}sz and Plummer.
  
  This version of the paper uses a different definition of a burl from
  the journal version of the paper~\cite{{esperetkkkn11}}, and a
  different proof of Lemma~\ref{l:3Burl}. This simplifies
  the exposition of our arguments throughout the whole
  paper.
\end{abstract}
\maketitle

\section{Introduction}

Given a graph $G$, let $\calM(G)$ denote the set of perfect matchings
in $G$.  A classical theorem of Petersen \cite{petersen1891} states
that every cubic bridgeless graph has at least one perfect matching,
i.e.\ $\calM(G)\neq \emptyset$.  Indeed, it can be proven that
any edge in a cubic bridgeless graph is contained in some perfect
matching~\cite{plesnik72}, which implies that $|\calM(G)|\geq 3$.

In the 1970s, Lov\'asz and Plummer conjectured that the number of
perfect matchings of a cubic bridgeless graph $G$ should grow
exponentially with its order (see~\cite[Conjecture
8.1.8]{lovaszp86}).  It is a simple exercise to prove that $G$
contains at most $2^{|V(G)|}$ perfect matchings, so we can state the
conjecture as follows:

\begin{lpc}
  There exists a universal constant $\epsilon >0$ such that for any
  cubic bridgeless graph $G$, $$2^{\epsilon |V(G)|} \leq |\calM(G)|
  \leq 2^{|V(G)|}.$$

\end{lpc}

The problem of computing $|\calM(G)|$ is connected to problems in
molecular chemistry and statistical physics (see e.g.~\cite[Section
8.7]{lovaszp86}). In general graphs, this problem is $\sharp
P$-complete \cite{valiant79}.  Thus we are interested in finding good
bounds on the number of perfect matchings for various classes of
graphs such as the bounds in the conjecture above.

For bipartite graphs, $|\calM(G)|$ is precisely the permanent of the
graph biadjacency matrix.  Voorhoeve proved the conjecture for cubic
bipartite graphs in 1979 \cite{voorhoeve79}; Schrijver later extended
this result to all regular bipartite graphs \cite{schrijver98}.  We
refer the reader to \cite{laurents10} for an exposition of this
connection and of an elegant proof of Gurvits generalizing Schrijver's
result. For {\em fullerene graphs}, a class of planar cubic graphs for
which the conjecture relates to molecular stability and aromaticity of
fullerene molecules, the problem was settled by Kardo\v{s}, Kr\'al',
Mi\v{s}kuf and Sereni \cite{kardoskms09}.  Chudnovsky and Seymour
recently proved the conjecture for all cubic bridgeless planar graphs
\cite{chudnovskys11}.

The general case has until now remained open.  Edmonds, Lov\'asz and
Pulleyblank \cite{edmondslp82} proved that any cubic bridgeless $G$
contains at least $\frac 14{|V(G)|}+2$ perfect matchings (see also
\cite{naddef82}); this bound was later improved to $\frac 12|V(G)|$
\cite{kralss09} and then $\frac34|V(G)|-10$ \cite{esperetkss10}.  The
order of the lower bound was not improved until Esperet, Kardo\v{s},
and Kr\'al' proved a superlinear bound in 2009 \cite{esperetkk11}.
The first bound, proved in 1982, is a direct consequence of a lower
bound on the dimension of the perfect matching polytope, while the
more recent bounds combine polyhedral arguments with analysis of
brick and brace decompositions.

In this paper we solve the general case.  To avoid technical
difficulties when contracting sets of vertices, we henceforth allow
graphs to have multiple edges, but not loops.  Let $m(G)$ denote
$|\calM(G)|$, and let $m^{\star}(G)$ denote the minimum, over all
edges $e\in E(G)$, of the number of perfect matchings containing $e$.
Our result is the following:

\begin{thm}\label{t:Main}
  For every cubic
  bridgeless graph $G$ we have $m(G) \geq 2^{|V(G)|/\ceps}$.
\end{thm}

We actually prove that at least one of two sufficient conditions applies:

\begin{thm}\label{t:Main2}
  For every cubic
  bridgeless graph $G$, at least one of the following holds:
\begin{itemize}
 \item[$\Sone$] $m^{\star}(G) \geq 2^{|V(G)|/\ceps},$ or
 \item[$\Stwo$] there exist $M,M' \in \calM(G)$ such that $M \triangle
   M'$ has at least $|V(G)|/\ceps$ components.
\end{itemize}
\end{thm}

To see that Theorem~\ref{t:Main2} implies Theorem~\ref{t:Main}, we can
clearly assume that $\Stwo$ holds since $m^{\star}(G) \leq m(G)$.
Choose $M,M' \in \calM(G)$ such that the set $\calC$ of components of
$M \triangle M'$ has cardinality at least $|V(G)|/\ceps$, and note that
each of these components is an even cycle alternating between $M$ and
$M'$.  Thus for any subset $\calC' \subseteq \calC$, we can construct
a perfect matching $M_{\calC'}$ from $M$ by flipping the edges on the
cycles in $\calC'$, i.e.\ $M_{\calC'} = M \triangle \bigcup_{C \in
  \calC'} C$. The $2^{|\calC|}$ perfect matchings $M_{\calC'}$ are
distinct, implying Theorem \ref{t:Main}.

\begin{figure}[t]
\centering
\includegraphics[width=8cm]{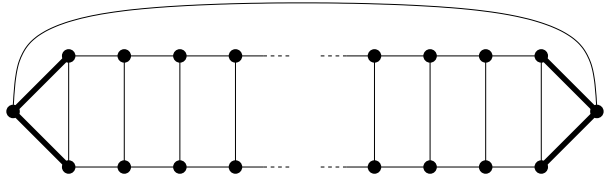}
\caption{A graph cubic bridgeless graph $G$ with
  $m^{\star}(G)=1$.} \label{f:uni}
\end{figure}

We cannot discard either of the sufficient conditions $\Sone$ or
$\Stwo$ in the statement of Theorem \ref{t:Main2}.  To see that
$\Stwo$ cannot be omitted, consider the graph depicted in
Figure~\ref{f:uni} and observe that each of the four bold edges is
contained in a unique perfect matching.  To see that $\Sone$ cannot be
omitted, it is enough to note that there exist cubic graphs with girth
logarithmic in their size (see~\cite{imrich84} for a construction).
Such graphs cannot have linearly many disjoint cycles, so condition
$\Stwo$ does not hold.

\subsection{Definitions and notation}
\

For a graph $G$ and a set $X \subseteq V(G)$, $G|X$ denotes the
subgraph of $G$ induced by $X$. For a set $X \subseteq V(G)$, let
$\delta(X)$ denote the set of edges with exactly one endpoint in $X$,
and let $E_X$ denote the set of edges with at least one endpoint in
$X$, i.e.\ $E_X = E(G|X)\cup \delta(X)$. The set $C=\delta(X)$ is
called an \emph{edge-cut}, or a \emph{$k$-edge-cut}, where $k=|C|$,
and $X$ and $V(G)\setminus X$ are the \emph{sides} of $C$.  A
$k$-edge-cut is said to be {\em even} (resp.\ {\em odd}) if $k$ is
even (resp.\ odd). Observe that the parity of an edge-cut $\delta(X)$
in a cubic graph is precisely that of $|X|$. An edge-cut $\delta(X)$
is \emph{cyclic} if both $G|X$ and $G|(V(G) \setminus X)$ contain a
cycle. Observe that every 2-edge-cut in a cubic graph is cyclic.  If
$G$ contains no edge-cut (resp.\ cyclic edge-cut) of size less than
$k$, we say that $G$ is \emph{$k$-edge-connected} (resp.\
\emph{cyclically $k$-edge-connected}).

Observe that the number of perfect matchings of a graph is the product
of the number of perfect matchings of its connected components. Hence,
in order to prove Theorem~\ref{t:Main}, we restrict ourselves to
connected graphs for the remainder of this paper (this means, for
example, that we can consider the terms \emph{2-edge-connected} and
\emph{bridgeless} to be interchangeable, and the sides of a cut are
well-defined).

\subsection{Constants}
\ 
Let $x:=\log(\tfrac43)/\log(2)$.
The following constants appear throughout the
paper: 
$$\alpha:=\tfrac{x}{314}, \qquad \beta_1 :=\tfrac{154x}{314},
\qquad \beta_2:=\tfrac{74x}{314}, \qquad \gamma:=\tfrac{312x}{314}.$$
We avoid using the numerical values of these constants for the sake of
clarity.  Throughout the paper we make use of the following
inequalities, which can be routinely verified:

\makeatletter
\let\temptag\tagform@
\def\tagform@#1{\maketag@@@{\bf(\ignorespaces#1\unskip\@@italiccorr)}}
\begin{align}
 0 < \alpha \leq \beta_2 \leq \beta_1,\label{e:const1}\\
 1/\ceps \leq \frac{\alpha}{9\beta_1+3},\label{e:const2}\\
 \beta_2+6\alpha \leq \beta_1,\label{e:const3}\\
 74\alpha \leq \beta_2,\label{e:const4}\\
 146 \alpha \leq \beta_1,\label{e:const5}\\
 \beta_2+ 80\alpha \leq \beta_1,\label{e:const6}\\
 6\alpha + \gamma \leq \log(6)/\log(2),\label{e:const7}\\
 \gamma + 2\beta_1 +7\alpha -\beta_2 \leq 1,\label{e:const8}\\
 6\alpha+2\beta_1 \leq \log(\tfrac{4}{3})/\log(2),\label{e:const9}\\
2\beta_1+4\alpha \leq \gamma.\label{e:const10}
\end{align}
\let\tagform@\temptag
\makeatother
The integer $\ceps$ is chosen minimum so that the system of inequalities
above has a solution.  Inequalities
{\bf (\ref{e:const4})}, {\bf (\ref{e:const6})}, {\bf (\ref{e:const9})}, and
{\bf (\ref{e:const10})} are tight.

\section{The proof of Theorem \ref{t:Main2}}

In this section we sketch the proof of Theorem \ref{t:Main2},
postponing the proofs of two main lemmas until later sections.  Our
general approach to Theorem \ref{t:Main2} is to reduce on cyclic
2-edge-cuts and cyclic 3-edge-cuts and prove inductively that either
$\Sone$ or $\Stwo$ holds.  Dealing with $\Sone$ is relatively
straightforward -- perfect matchings containing a given edge behave
well with reductions on a cut, which is our main motivation for
considering $m^{\star}(G)$.  To deal with $\Stwo$, we do not directly
construct perfect matchings $M$ and $M'$ for which $M\triangle M'$ has
many components.  Instead, we define a special type of vertex set in
which a given random perfect matching is very likely to admit an
alternating cycle.  We call these sets {\em burls} and we call a set
of disjoint burls a {\em foliage} -- a large foliage will guarantee
the existence of two perfect matchings with many components in their
symmetric difference.

\subsection{Burls, twigs, and foliage weight}
\

Consider a subset $X \subseteq V(G)$. Let $\calM(G,X)$ denote the
family of subsets $M$ of $E_X$ (the edges with at least one endpoint
in $X$) such that every vertex of $X$ is incident with exactly one
edge of $M$. Note that some elements of $\calM(G,X)$ might not be
matchings in $G$ (if two edges of $\delta(X)$ share a vertex from
$V(G)\setminus X$). However, for any $M \in \calM(G,X)$, $M \cap
E(G|X)$ is a matching. 

A probability distribution $\bM$ on $\calM(G,X)$ is \emph{balanced} if
for any edge $e \in E_X$, $\Pr[e\in \bM] = \tfrac13$. It follows from
Edmonds' characterization of the perfect matching
polytope~\cite{edmonds65} that if $G$ is cubic and bridgeless, there
exists a balanced probability distribution on
$\calM(G,V(G))=\calM(G)$. For any $X \subseteq V(G)$, the restriction
of this distribution to $E_X$ yields a balanced probability
distribution on $\calM(G,X)$.  The following easy fact will be used
several times throughout the proof:

\begin{cla}\label{cl:balanced}
  Let $G$ be a cubic bridgeless graph and consider $Y \subseteq X
  \subseteq V(G)$ such that $C=\delta(Y)$ is a 3-edge-cut in $G$. For
  any balanced probability distribution $\bM$ on $\calM(G,X)$, and any
  $M\in \calM(G,X)$ such that $\Pr[\bM =M]>0$, we have $|M\cap C|=1$.
\end{cla}

Given some $M \in \calM(G,X)$, a cycle of $G|X$ is
\emph{$M$-alternating} if it has even length and half of its edges are
in $M$ (it alternates edges in $M$ and edges not in $M$). Let
$a(G,X,M)$ denote the maximum number of disjoint $M$-alternating
cycles in $G|X$ (equivalently, the maximum number of components of $M
\triangle M'$, for $M' \in \calM(G,X)$).

We define a \emph{burl} as a vertex set $X \subseteq V(G)$ such that
for any balanced probability distribution $\bM$ on $\calM(G,X)$,
$\Ex[a(G,X,\bM)] \ge \tfrac13$. Note that if $X$ is a burl, any set $Y
\supset X$ is also a burl, since any balanced probability distribution
on $\calM(G,Y)$ induces a balanced probability distribution on
$\calM(G,X)$. We would like to insist on the fact that we consider the
whole set $\calM(G,X)$, and not only $\{M \cap E_X, M \in
\calM(G)\}$. This way, being a burl is really a local property of $X$
and is completely independent of the structure of $G|(V(G) \setminus
X)$. This aspect of burls will be fundamental in the proof of
Theorem~\ref{t:Main2}.

A collection of disjoint vertex sets $\{X_1,\ldots,X_k\}$ is a {\em
  foliage} if each $X_i$ is a burl.  Assume that $G$ contains such a
collection of disjoint sets, and consider a balanced probability
distribution $\bM$ on $\calM(G,V(G))=\calM(G)$. This distribution
induces balanced probability distributions $\bM_{X_i}$ on
$\calM(G,X_i)$, for each $1\le i \le k$. By definition of a burl, we
have $\Ex[a(G,X_i,\bM_{X_i})] \ge \tfrac13$ for each each $1\le i \le
k$. By linearity of expectation, the maximum number of disjoint
alternating cycles of $\bM$ is then expected to be at least $k/3$. We
get the following key fact as a consequence:

\begin{cor}\label{c:FoliageApplied}
  If a cubic bridgeless graph $G$ contains a foliage $\calX$, then
  there exist perfect matchings $M,M' \in \calM(G)$ such that $M
  \triangle M'$ has at least $|\calX|/3$ components.
\end{cor}

We now introduce a special class of burls. Let $G$ be a cubic
bridgeless graph and let $X \subseteq V(G)$. We say that $X$ is a
\emph{$2$-twig} if $|\delta(X)|=2$, and $X$ is a \emph{$3$-twig} if
$|\delta(X)|=3$ and $|X| \geq 5$ (that is, $X$ is neither a triangle,
nor a single vertex).  A \emph{twig} in $G$ is a $2$- or $3$-twig.
Before we prove that every twig is a burl, we need a simple lemma.

\begin{lem}\label{l:smallcount}
Let $G$ be a cubic bridgeless graph. Then
\begin{enumerate}
 \item $m(G-e) \geq 2$ for every $e \in E(G)$, and
 \item $m(G) \geq 4$ if $|V(G)|\geq 6$.  In particular, for any $v\in
   V(G)$ there is an $e\in \delta(\{v\})$ contained in at least two
   perfect matchings.
\end{enumerate}
\end{lem}
\begin{proof}
  The first item follows from the classical result mentioned in the
  introduction: every edge of a cubic bridgeless graph is contained in
  a perfect matching.  The second is implied by the bound $m(G)\geq
  \tfrac 14|V(G)|+2$ from~\cite{edmondslp82}.
\end{proof}

\begin{lem}\label{l:twigs}
Every twig $X$ in a cubic bridgeless graph $G$ is a burl.
\end{lem}
\begin{proof}
  Let $\bM$ be a balanced probability distribution on $\calM(G,X)$.

  If $X$ is a $2$-twig, let $H$ be obtained from $G|X$ by adding an
  edge $e$ joining the two vertices incident with $\delta(X)$. Then
  $H$ is cubic and bridgeless. By applying Lemma~\ref{l:smallcount}(1)
  to $H$, we see that $G|X$ contains at least one $M$-alternating
  cycle for every $M \in \calM(G,X)$ such that $M \cap
  \delta(X)=\emptyset$. Note that since $H$ is cubic, $|X|$ is even,
  and thus $M$ either contains the two edges of $\delta(X)$, or none
  of them. Since $\bM$ is balanced, $\Pr[\bM \cap \delta(X)=
  \emptyset] \geq 1-1/3=2/3$. Hence $\Ex[a(G,X,\bM)] \ge \tfrac23$ and
  we conclude that $X$ is a burl.

  Suppose now that $X$ is a $3$-twig. Let
  $\delta(X)=\{e_1,e_2,e_3\}$. Let $H$ be obtained from $G$ by
  identifying all the vertices in $V(G)-X$ (removing loops but
  preserving multiple edges). We apply Lemma~\ref{l:smallcount}(2) to
  $H$, which is again cubic and bridgeless. It follows that for some
  $1 \leq i \leq 3$, the edge $e_i$ is in at least two perfect
  matchings of $H$. Therefore $G|X$ contains at least one
  $M$-alternating cycle for every $M \in \calM(G,X)$ such that $M \cap
  \delta(X)= \{e_i\}$. By Claim~\ref{cl:balanced}, $\Pr[\bM \cap
  \delta(X)= \{e_i\}]=\Pr[e_i \in \bM]=1/3$. It implies that
  $\Ex[a(G,X,\bM)] \ge \tfrac13$ and thus $X$ is a burl.
\end{proof}

The \emph{weight} of a foliage $\calX$ containing $k$ twigs is defined
as $\fw(\calX) := \beta_1 k + \beta_2 (|\calX|-k)$, that is each twig
has weight $\beta_1$ and each non-twig burl has weight $\beta_2$. Let
$\fw(G)$ denote the maximum weight of a foliage in a graph $G$.

\subsection{Reducing on small edge-cuts}
\

We now describe how we reduce on 2-edge-cuts and 3-edge-cuts, and
consider how these operations affect $m^{\star}(G)$ and foliages.  Let
$C$ be a $3$-edge-cut in a cubic bridgeless graph $G$. The two graphs
$G_1$ and $G_2$ obtained from $G$ by identifying all vertices on one
of the sides of the edge-cut (removing loops but preserving multiple
edges) are referred to as \emph{$C$-contractions} of $G$ and the
vertices in $G_1$ and $G_2$ created by this identification are called
\emph{new}.

We need a similar definition for $2$-edge-cuts. Let $C=\{e,e'\}$ be a
$2$-edge-cut in a cubic bridgeless graph $G$. The two
\emph{$C$-contractions} $G_1$ and $G_2$ are now obtained from $G$ by
deleting all vertices on one of the sides of $C$ and adding an edge joining
the remaining ends of $e$ and $e'$. The resulting edge is now called
\emph{new}.

In both cases we say that $G_1$ and $G_2$ \emph{are obtained from $G$
  by a cut-contraction}. The next lemma provides some useful
properties of cut-contractions.

\begin{lem}\label{l:Contraction}
  Let $G$ be a graph, let $C$ be a $2$- or a $3$-edge-cut in $G$, and
  let $G_1$ and $G_2$ be the two $C$-contractions. Then
\begin{enumerate}
\item $G_1$ and $G_2$ are cubic bridgeless graphs,
  \item $m^{\star}(G) \geq m^{\star}(G_1)\,m^{\star}(G_2)$, and
  \item For $i=1,2$ let $\calX_i$ be a foliage in $G_i$ such that for
    every $X \in \calX_i$, if $|C|=3$ then $X$ does not contain the
    new vertex, and if $|C|=2$ then $E(G_i|X)$ does not contain the
    new edge.  Then $\calX_1 \cup \calX_2$ is a foliage in $G$. In
    particular, we have $\fw(G) \geq \fw(G_1)+\fw(G_2)-2\beta_1$.
\end{enumerate}
\end{lem}
\begin{proof}\mbox{}\smallskip
\begin{enumerate}

\item This can be confirmed routinely.

\item Consider first the case of the contraction of a 2-edge-cut
  $C=\delta(X)$ in $G$. Let $e$ be an edge with both ends in
  $X=V(G_1)$. Every perfect matching of $G_1$ containing $e$ combines
  either with $m^{\star}(G_2)$ perfect matchings of $G_2$ containing
  the new edge of $G_2$, or with $2m^{\star}(G_2)$ perfect matchings
  of $G_2$ avoiding the new edge of $G_2$. If $e$ lies in $C$, note
  that perfect matchings of $G_1$ and $G_2$ containing the new edges
  can be combined into perfect matchings of $G$ containing $C$. Hence,
  $e$ is in at least $m^{\star}(G_1)\,m^{\star}(G_2)$ perfect matchings
  of $G$.

  Now consider a 3-edge-cut $C=\delta(X)$. If $e$ has both ends in $X
  \subset V(G_1)$, perfect matchings of $G_1$ containing $e$ combine
  with perfect matchings of $G_2$ containing either of the 3 edges of
  $C$. If $e$ is in $C$, perfect matchings containing $e$ in $G_1$ and
  $G_2$ can also be combined into perfect matchings of $G$. In any
  case, $e$ is in at least $m^{\star}(G_1)\,m^{\star}(G_2)$ perfect
  matchings of $G$.

\item In this case the elements of $\calX_1 \cup \calX_2$ are disjoint
  subsets of $V(G)$. Consider some $X \in \calX_1 \cup \calX_2$, and
  assume without loss of generality that $X \in \calX_1$. Since $X$
  does not contain the new vertex of $G_1$ (if $|C|=3$), or the new
  edge (if $|C|=2$), each balanced probability distribution on
  $\calM(G,X)$ is also a balanced probability distribution on
  $\calM(G_1,X)$, so $\calX_1 \cup \calX_2$ is a foliage in $G$. Since
  $\beta_1 \geq \beta_2$, this implies $\fw(G) \geq
  \fw(G_1)+\fw(G_2)-2\beta_1$.  \qedhere
\end{enumerate}
\end{proof}

It is not generally advantageous to reduce on a $3$-edge-cut arising
from a triangle, unless this reduction leads to a chain of similar
reductions.  Thus we wish to get rid of certain triangles from the
outset.  We say that a triangle sharing precisely one edge with a
cycle of length three or four in a graph $G$ is \emph{relevant}, and
otherwise it is \emph{irrelevant}.  A graph $G$ is \emph{pruned} if it
contains no irrelevant triangles.  The following easy lemma shows that
we can prune a bridgeless cubic graph by repeated cut-contraction
without losing too many vertices.

\begin{lem}\label{l:ContractTriangles}
  Let $G$ be a cubic bridgeless graph, and let $k$ be the size of
  maximum collection of vertex-disjoint irrelevant triangles in
  $G$. Then one can obtain a pruned cubic bridgeless graph $G'$ from
  $G$ with $|V(G')| \geq |V(G)|-2k$ by repeatedly contracting
  irrelevant triangles.
\end{lem}

\begin{proof}
  We proceed by induction on $k$. Let a graph $G''$ be obtained from
  $G$ by contracting an irrelevant triangle $T$. The graph $G''$ is
  cubic and bridgeless by Lemma~\ref{l:Contraction}(1).  Since $T$ is
  irrelevant in $G$, the unique vertex of $G''$ obtained by
  contracting $T$ is not in a triangle in $G''$. Therefore if
  $\mathcal{T}$ is a collection of vertex disjoint irrelevant
  triangles in $G''$ then $\mathcal{T} \cup \{T\}$ is such a
  collection in $G$. (After the contraction of an irrelevant triangle,
  triangles that were previously irrelevant might become relevant, but
  the converse is not possible.) It follows that $|\mathcal{T}| \leq
  k-1$. By applying the induction hypothesis to $G''$, we see that the
  lemma holds for $G$.
\end{proof}

\begin{cor}\label{c:Prune}
  Let $G$ be a cubic bridgeless graph. Then we can obtain a cubic
  bridgeless pruned graph $G'$ from $G$ with $|V(G')| \geq |V(G)|/3$
  by repeatedly contracting irrelevant triangles.
\end{cor}

We wish to restrict our attention to pruned graphs, so we must make
sure that the function $m^{\star}(G)$ and the maximum size of a
foliage does not increase when we contract a triangle.

\begin{lem}\label{l:Triangle}
  Let $G'$ be obtained from a graph $G$ by contracting a
  triangle. Then $m^{\star}(G') \leq m^{\star}(G)$ and the maximum
  size of a foliage in $G'$ is at most the maximum size of a foliage
  in $G$.
\end{lem}
\begin{proof}
  Let $xyz$ be the contracted triangle, and let $e_x$, $e_y$, and
  $e_z$ be the edges incident with $x$, $y$, $z$ and not contained in
  the triangle in $G$. Let $t$ be the vertex of $G'$ corresponding to
  the contraction of $xyz$. Every perfect matching $M'$ of $G'$ has a
  canonical extension $M$ in $G$: assume without loss of generality
  that $e_x$ is the unique edge of $M'$ incident to $t$.  Then $M$
  consists of the union of $M'$ and $yz$. Observe that perfect
  matchings in $G$ containing $yz$ necessarily contain $e_x$, so every
  edge of $G$ is contained in at least $m^{\star}(G')$ perfect
  matchings.

  Now consider a burl $X'$ in $G'$ containing $t$. We show that $X=X'
  \cup \{x,y,z\} \setminus t$ is a burl in $G$. Let $\bM$ be a
  balanced probability distribution on $\calM(G,X)$. By
  Claim~\ref{cl:balanced} and the remark above, we can associate a
  balanced probability distribution $\bM'$ on $\calM(G',X')$ to $\bM$
  such that $\Ex[a(G,X,\bM)] =\Ex[a(G',X',\bM')]$. Since $X'$ is a
  burl in $G'$, this expectation is at least $\tfrac13$ and $X$ is a
  burl in $G$.

Since a burl avoiding $t$ in $G'$ is also a burl in $G$, it follows
  from the analysis above that the maximum size of a foliage cannot
  increase when we contract a triangle.
\end{proof}

\subsection{Proving Theorem \ref{t:Main2}}
\

We say that $G$ \emph{has a core} if we can obtain a cyclically
4-edge-connected graph $G'$ with $|V(G')|\geq 6$ by applying a
(possibly empty) sequence of cut-contractions to $G$ (recall that this
notion was defined in the previous subsection).

We will deduce Theorem \ref{t:Main2} from the next two lemmas.  This
essentially splits the proof into two cases based on whether or not
$G$ has a core.

\begin{lem}\label{l:Klee}
  Let $G$ be a pruned cubic bridgeless graph. Let $Z \subseteq V(G)$
  be such that $|Z| \geq 2$ and $|\delta(Z)| = 2$, or $|Z| \geq 4$ and
  $|\delta(Z)| = 3$. Suppose that the $\delta(Z)$-contraction $G'$ of
  $G$ with $Z \subseteq V(G')$ has no core. Then there exists a
  foliage $\calX$ in $G$ with $\bigcup_{X \in \calX}X \subseteq Z$
  and $$\fw(\calX) \geq \alpha|Z| + \beta_2.$$
\end{lem}

By applying Lemma~\ref{l:Klee} to a cubic graph $G$ without a core and
$Z=V(G) \setminus \{v\}$ for some $v \in V(G)$, we obtain the following.

\begin{cor}\label{c:Klee}
  Let $G$ be a pruned cubic bridgeless graph without a core.
  Then $$\fw(G) \geq \alpha (|V(G)|-1) + \beta_2.$$
\end{cor}

On the other hand, if $G$ has a core, we will prove that either $\fw(G)$
is linear in the size of $G$ or every edge of $G$ is contained in an
exponential number of perfect matchings.

\begin{lem}\label{l:MainInduction} Let $G$ be a pruned cubic
  bridgeless graph. If $G$ has a core then
$$m^{\star}(G) \geq 2^{\alpha|V(G)|-\fw(G)+\gamma}.$$
\end{lem}

We finish this section by deriving Theorem~\ref{t:Main2} from
Lemmas~\ref{l:Klee} and~\ref{l:MainInduction}. 
\begin{proof}[Proof of Theorem~\ref{t:Main2}] Let $\epsilon
  :=1/\ceps$.  By Corollary~\ref{c:Prune} there exists a pruned cubic
  bridgeless graph $G'$ with $|V(G')| \geq |V(G)|/3$ obtained from $G$
  by repeatedly contracting irrelevant triangles. Suppose first that
  $G'$ has a core. By Corollary~\ref{c:Prune} and Lemmas~\ref{l:Triangle}
  and~\ref{l:MainInduction}, condition $\Sone$ holds as long as
  $\epsilon|V(G)| \leq \alpha|V(G)|/3-\fw(G')$. Therefore we assume
  $\fw(G') \geq (\tfrac \alpha 3-\epsilon)|V(G)|$. It follows from the
  definition of $\fw(G')$ that $G'$ has a foliage containing at least
  $(\tfrac \alpha 3-\epsilon)|V(G)|/\beta_1$ burls. If $G'$ has no
  core then by Corollary~\ref{c:Klee} and the fact that $\alpha \leq
  \beta_2$, $\fw(G')\ge \alpha (|V(G')|-1)+\beta_2 \ge \alpha|V(G')|$,
  so $G'$ contains a foliage of size at least $\alpha|V(G')|/\beta_1
  \ge \alpha|V(G)|/3\beta_1$.  In both cases condition $\Stwo$ holds
  by Corollary~\ref{c:FoliageApplied} and Lemma~\ref{l:Triangle},
  since Equation {\bf (\ref{e:const2})} tells us that $3\epsilon \leq
  {(\tfrac \alpha3 -\epsilon)}/\beta_1$.
\end{proof}

\section{Cut decompositions}\label{s:CutDecompose}

In this section we study cut decompositions of cubic bridgeless
graphs. We mostly follow notation from~\cite{chudnovskys11}, however we
consider $2$- and $3$-edge-cuts simultaneously. Cut decompositions play a
crucial role in the proof of Lemma~\ref{l:Klee} in the next section.

Let $G$ be a graph. A \emph{non-trivial cut-decomposition} of $G$ is a
pair $(T, \phi)$ such that:
\begin{itemize}
 \item $T$ is a tree with $E(T) \neq \emptyset$,
 \item $\phi : V(G) \to V (T)$ is a map, and
 \item $|\phi^{-1}(t)|+\deg_T(t) \geq 3$ for each $t \in V(T)$. 
\end{itemize}

For an edge $f$ of $T$, let $T_1$, $T_2$ be the two components of $T
\setminus f$, and for $i = 1,2$ let $X_i = \phi^{-1}(T_i)$. Thus
$(X_1,X_2)$ is a partition of $V(G)$ that induces an edge-cut
denoted by $\phi^{-1}(f)$. If $|\phi^{-1}(f)| \in \{2,3\}$ for each $f
\in E(T)$ we call $(T, \phi)$ a \emph{small-cut-decomposition} of $G$.

Let $(T, \phi)$ be a small-cut-decomposition of a $2$-edge-connected
cubic graph $G$, and let $T_0$ be a subtree of $T$ such that
$\phi^{-1}(V(T_0)) \neq \emptyset$. Let $T_1,\ldots,T_s$ be the
components of $T \setminus V (T_0)$, and for $1 \leq i \leq s$ let
$f_i$ be the unique edge of $T$ with an end in $V(T_0)$ and an end in
$V(T_i)$. For $0 \leq i \leq s$, let $X_i = \phi^{-1}(V(T_i))$. Thus
$X_0,X_1, \ldots ,X_s$ form a partition of $V(G)$.
Let $G'$ be the graph obtained from $G$ as follows.  Set $G_0=G$.  For
$i=1,\ldots,s$, take $G_{i-1}$ and let $G_i$ be the
$(\phi^{-1}(f_i))$-contraction containing $X_0$. Now let $G'$ denote
$G_{s}$.
Note that $G'$ is cubic. We call $G'$ \emph{the hub of $G$ at $T_0$}
(with respect to $(T, \phi)$). If $t_0 \in V (T)$ and $\phi^{-1}(t_0)
\ne \emptyset$, by the \emph{hub of $G$ at $t_0$} we mean the hub of
$G$ at $T_0$, where $T_0$ is the subtree of $T$ with vertex set
$\{t_0\}$.

Let $\calY$ be a collection of disjoint subsets of $V(G)$. We say that
a small-cut-decomposition $(T,\phi)$ of $G$ \emph{refines} $\calY$ if
for every $Y \in \calY$ there exists a leaf $v \in V(T)$ such that $Y
= \phi^{-1}(v)$. Collections of subsets of $V(G)$ that can be refined
by a small-cut decomposition are charaterized in the following easy
lemma.

\begin{lem}\label{l:RefinableCollections} Let $G$ be a cubic
  bridgeless graph.  Let $\calY$ be a collection of disjoint subsets
  of $V(G)$. Then there exists a small-cut-decomposition refining
  $\calY$ if $|Y| \geq 2$ and $|\delta(Y)| \in \{2,3\}$ for every $Y
  \in \calY$, and either
\begin{enumerate}
 \item  $\calY =\emptyset$ and $G$ is not cyclically $4$-edge-connected, or
 \item  $\calY=\{Y\}$, and $|V(G) \setminus  Y| >1$, or
 \item  $|\calY| \geq 2$. 
\end{enumerate}
\end{lem}
\begin{proof}
  We only consider the case $|\calY|\geq 3$, as the other cases are
  routine. Take $T$ to be a tree on $|\calY|+1$ vertices with
  $|\calY|$ leaves $\{v_Y \: | \: Y \in \calY \}$ and a non-leaf
  vertex $v_0$. The map $\phi$ is defined by $\phi(u)=v_Y$, if $u \in
  Y$ for some $Y \in \calY$, and $\phi(u)=v_0$, otherwise.  Clearly,
  $(T,\phi)$ refines $\calY$ and is a small-cut-decomposition of $G$.
\end{proof}
 
We say that $(T,\phi)$ is \emph{$\calY$-maximum} if it refines $\calY$
and $|V(T)|$ is maximum among all small-cut decompositions of $G$
refining $\calY$. The following lemma describes the structure of
$\calY$-maximum decompositions. It is a variation of Lemma 4.1 and
Claim 1 of Lemma 5.3 in~\cite{chudnovskys11}.

\begin{lem}\label{l:CutDecompose} Let $G$ be a cubic bridgeless
  graph. Let $\calY$ be a collection of disjoint subsets of $V(G)$ and
  let $(T,\phi)$ be a $\calY$-maximum small-cut-decomposition of
  $G$. Then for every $t \in V(T)$ either $\phi^{-1}(t)=\emptyset$, or
  $\phi^{-1}(t)\in \calY$, or the hub of $G$ at $t$ is cyclically
  $4$-edge-connected.
\end{lem}

\begin{proof}
  Fix $t \in V(T)$ with $\phi^{-1}(t) \neq \emptyset$ and
  $\phi^{-1}(t) \not \in \calY$. Let $f_1,\ldots, f_k$ be the edges
  of $T$ incident with $t$, and let $T_1,\ldots, T_k$ be the
  components of $T \setminus \{t\}$, where $f_i$ is incident with a
  vertex $t_i$ of $T_i$ for $1 \leq i \leq k$. Let $X_0$ =
  $\phi^{-1}(t)$, and for $1 \leq i \leq k$ let $X_i =
  \phi^{-1}(V(T_i))$. Let $G'$ be the hub of $G$ at $t$, and let $G''$
  be the graph obtained from $G'$ by subdividing precisely once every
  new edge $e$ corresponding to the cut-contraction of a cut $C$ with
  $|C|=2$. The vertex on the subdivided edge $e$ is called \emph{the
    new vertex corresponding to the cut-contraction of $C$}, by
  analogy with the new vertex corresponding to the cut-contraction of
  a cyclic 3-edge-cut.

  Note that $G'$ is cyclically 4-edge-connected if and only if $G''$
  is cyclically 4-edge-connected.  Suppose for the sake of
  contradiction that $C=\delta(Z)$ is a cyclic edge-cut in $G''$ with
  $|C|\leq 3$. Then $|C|\in \{2,3\}$ by Lemma~\ref{l:Contraction}(1),
  as $G''$ is a subdivision of $G'$ and $G'$ can be obtained from $G$
  by repeated cut-contractions.  Let $T'$ be obtained from $T$ by by
  splitting $t$ into two vertices $t'$ and $t''$, so that $t_i$ is
  incident to $t'$ if and only if the new vertex of $G''$
  corresponding to the cut-contraction of $\phi^{-1}(f_i)$ is in
  $Z$. Let $\phi'(t')=X_0 \cap Z$, $\phi'(t'')=X_0  \setminus  Z$, and
  $\phi'(s)=\phi(s)$ for every $s \in V(T') \setminus \{t',t''\}$.

  We claim that $(T',\phi')$ is a small-cut-decomposition of $G$
  contradicting the choice of $T$. It is only necessary to verify that
  $|\phi^{-1}(s)|+\deg_{T'}(s)\geq 3$ for $s \in \{t',t''\}$. We have
  $|\phi^{-1}(t')|+\deg_{T'}(t')-1=|Z \cap V(G'')| \geq 2$ as $C$ is a
  cyclic edge-cut in $G''$. It follows that
  $|\phi^{-1}(t')|+\deg_{T'}(t') \geq 3$ and the same holds for $t''$
  by symmetry.
\end{proof}

\begin{figure}[htbp]
\centering
\includegraphics[scale=0.5]{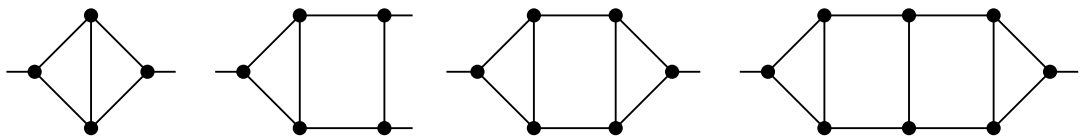}
\caption{Isomorphism classes of subgraphs induced by elementary
  twigs.} \label{f:twigs}
\end{figure}

We finish this section by describing a collection $\calY$ to which we
will be applying Lemma~\ref{l:CutDecompose} in the sequel.  In a cubic
bridgeless graph $G$ a union of the vertex set of a relevant triangle
with the vertex set of a cycle of length at most four sharing an edge
with it is called \emph{a simple twig}. Note that simple twigs
corresponding to distinct relevant triangles can intersect, but one
can routinely verify that each simple twig intersects a simple twig
corresponding to at most one other relevant triangle. \emph{An
  elementary twig} is either a simple twig, that intersects no simple
twig corresponding to a relevant triangle not contained in it, or the
union of two intersecting simple twigs, corresponding to distinct
relevant triangles. An elementary twig is, indeed, a twig, unless it
constitutes the vertex set of the entire graph. Figure~\ref{f:twigs}
shows all possible elementary twigs. The next corollary follows
immediately from the observations above and
Lemmas~\ref{l:RefinableCollections} and~\ref{l:CutDecompose}.

\begin{cor}\label{c:ElemTwigs}
 Let $G$ be a cubic bridgeless graph that is not cyclically
 $4$-edge-connected with $|V(G)| \geq 8$. Then there exists a
collection $\calY$ of pairwise disjoint
 elementary twigs in $G$ such that every relevant triangle in $G$ is
contained in an element of $\calY$.
 Further, there exists a $\calY$-maximum
 small-cut-decomposition $(T, \phi)$ of $G$ and for every $t \in
 V(T)$ either $\phi^{-1}(t)=\emptyset$, or $\phi^{-1}(t)$ is an
 elementary twig, or the hub of $G$ at $t$ is cyclically
 $4$-edge-connected.
\end{cor}

\section{Proof  of Lemma~\ref{l:Klee}.}\label{s:Klee}

The proof of Lemma~\ref{l:Klee} is based on our ability to find burls
locally in the graph. The following lemma is a typical example.

\begin{lem}\label{l:4Burl}
  Let $G$ be a cubic bridgeless graph and let $X \subseteq V(G)$ be
  such that $|\delta(X)|=4$ and $m(G|X) \geq 2$. Then $X$ is a burl.
\end{lem}
\begin{proof}
  Note that if $M \in \calM(G,X)$ contains no edges of $\delta(X)$
  then $G|X$ contains an $M$-alternating cycle. Let $\bM$ be a balanced
  probability distribution on $\calM(G,X)$. As $M \cap \delta(X)$ is
  even for every $M \in \calM(G,X)$ we have $$\tfrac 43=\Ex\left[|\bM
    \cap \delta(X)|\right] \geq 2\Pr[\bM \cap \delta(X) \neq
  \emptyset].$$ Therefore $\Pr[\bM \cap \delta(X) = \emptyset] \geq
  1/3$. Hence, $\Ex[a(G,X,\bM)] \ge \tfrac13$ and so $X$ is a burl.
\end{proof}

The proof of Lemma~\ref{l:Klee} relies on a precise study of the
structure of small-cut trees for graphs with no core. The following
two lemmas indicate that long paths in such trees necessarily contain
some burls.

\begin{lem}\label{l:3Burl} Let $(T,\phi)$ be a small-cut-decomposition
  of a cubic bridgeless graph $G$, and let $P$ be a path in $T$ with
  $|V(P)|= 10$. If we have
 \begin{itemize}
 \item $\deg_T(t)=2$ for every $t \in V(P)$,
 \item the hub of $G$ at $t$ is isomorphic to $K_4$ for every $t \in
   V(P)$, and
   \item $|\phi^{-1}(f)|=3$ for every edge $f \in E(T)$ incident to a
     vertex in $V(P)$,
 \end{itemize}
then $\phi^{-1}(P)$ is a burl.
\end{lem}

\begin{proof}
  Let $P'=v_{-1}v_0\ldots v_{9}v_{10}$ be a path in $T$ such that
  $P=v_0\ldots v_{9}$.  Let $f_i=v_{i-1}v_{i}$ and let
  $C_i=\{e^i_1,e^i_2,e^i_3\}=\phi^{-1}(f_i)$, $0 \le i \le 10$. Let
  $X:=\phi^{-1}(V(P))$.  It is easy to see that $\phi^{-1}(v_i)$
  contains precisely two vertices joined by an edge, $0\le i \le 9$.

We assume without loss of generality that $G|X$ contains no cycles of
length $4$, as otherwise the lemma holds by Lemma~\ref{l:4Burl}.  Let
$A$ be the set of ends of edges in $C_{0}$ outside of $X$, and let $B$
be the set of ends of edges in $C_{10}$ outside of $X$.  Observe that
$E_X$ consists of $3$ internally vertex-disjoint paths from $A$ to
$B$, as well as one edge in $G|\phi^{-1}(\{v_i\})$ for each $0\leq i\leq
9$.  Let $R_1$, $R_2$ and $R_3$ be these three paths from $A$ to $B$,
and let $u_j$ and $v_j$ be the ends of $R_j$ in $A$ and $B$,
respectively, for $j=1,2,3$.  For $0 \leq i \leq 9$, we have
$\phi^{-1}(v_i)=\{x_i,y_i\}$ so that $x_i \in V(R_j), y_i \in
V(R_{j'})$ for some $\{j,j'\} \subseteq \{1,2,3\}$ with $j\ne j'$
  ; let the \emph{index} $\sigma_i$ of $v_i$ be defined as
  $\{j,j'\}$. Since there is no 4-cycle in $X$, $\sigma_i\ne
  \sigma_{i-1}$ for $1\le i\le 9$.  Let the \emph{type} $\psi_i$ of
  $v_i$ (for $1\le i \le 8$) be defined as 0 if
  $\sigma_{i-1}=\sigma_{i+1}$, otherwise let $\psi_i=1$.
  
  Let $i,j,k$ be integers such that $0\le i < j< k\le 10$
  which will be determined later. Let
  $X_1:=\phi^{-1}(\{v_i,\dots,v_{j-1}\})$,
  $X_2:=\phi^{-1}(\{v_j,\dots,v_{k-1}\})$, $X_0=X_1\cup X_2$.

  Let $\bM$ be a balanced probability distribution on $\calM(G,X)$,
  let $Z_0$ ($Z_1$, $Z_2$) be the maximum number of disjoint
  $\bM$-alternating cycles in $G|X_0$ ($G|X_1$, $G|X_2$,
  respectively). Let $Y_\ell=|\bM\cap C_\ell|$, for every $\ell$, and
  let $Y=Y_i+Y_j+Y_k$. Since $\bM$ is balanced, we have $\Ex(Y)=3$;
  moreover, $Y_i\equiv Y_j \equiv Y_k \pmod{2}$. Therefore,
  $\Pr(Y=1)=0$; $Y=3$ if and only if $Y_i=Y_j=Y_k=1$; and $ Y=2$ if
  and only if $\{Y_i,Y_j,Y_k\}=\{2,0,0\}$.

Assume that $i,j,k$ fulfill the following conditions:
\begin{enumerate}
\item $\Pr(Z_1=0\,|\,Y_i=Y_j=0)=0$, $\Pr(Z_2=0\,|\,Y_j=Y_k=0)=0$, and
  $\Pr(Z_0=0\,|\,Y_i=Y_k=0)=0$;
\item for at least one of the cuts $C_i$, $C_j$ or $C_k$, say $C_t$,
  there exists an edge $e\in C_t$ such that for at least one of the
  two corresponding graphs among $G|X_0$, $G|X_1$,
  $G|X_2$, say $G|X_s$, there is an alternating cycle in
  $G|X_s$ for any element of $\calM(G,X)$ containing $e$, provided
  $Y_i=Y_j=Y_k=1$.
\end{enumerate}

First, we derive $\Ex(Z_0)\ge \frac13$ from these assumptions, then we
prove the existence of such a triple $i,j,k$.

Observe that the first condition yields $\Ex(Z_0\,|\,Y=0)\ge 2$ and
$\Ex(Z_0\,|\,Y=2)\ge 1$. Since $\Ex(Y)=3$, we have $3\cdot
\Pr(Y=0)+\Pr(Y=2) \ge \Pr(Y\ge 4)$. This gives $\Pr(Y\ne 3)\le 4\cdot
\Pr(Y=0) + 2\cdot \Pr(Y=2)$, and hence $\Ex(Z_0\,|\,Y\ne 3)\ge
\frac12$. Let $C_t=\{e_1,e_2,e_3\}$, where $e=e_1$. Let
$p_i=\Pr[\bM\cap C_t=\{e_i\} \wedge Y=3]$, $i=1,2,3$. Clearly
$p_1+p_2+p_3=\Pr(Y=3)$. On the other hand, since $\bM$ is balanced,
$\frac13-p_1 \le \frac13-p_2+\frac13-p_3$ (all elements of
$\calM(G,X)$ containing $e_1$ together with some other edge from $C_t$
contain $e_2$ or $e_3$). Hence, $p_1\ge \frac12\cdot
(\Pr(Y=3)-\frac13)$.

Altogether, in this case
$$
\gathered \Ex(Z_0)=\Ex(Z_0\,|\,Y\ne 3)\cdot \Pr(Y\ne
3)+\Ex(Z_0\,|\,Y=3)\cdot \Pr(Y=3)\ge \\\ge \tfrac12\cdot
\left(1-\Pr(Y=3)\right)+\tfrac12\cdot \left(\Pr(Y=3)-\tfrac13\right) =
\tfrac13.
\endgathered
$$

Now we prove that there is always $i,j,k$ such that both conditions
(1) and (2) are satisfied. Note that (1) is satisfied if $j-i=3$ and
$\psi_{i+1}=1$, or if $j-i\ge 4$; the same holds for $j$ and $k$ (to
observe this, see Figure \ref{fig:tbd}).

\begin{figure}
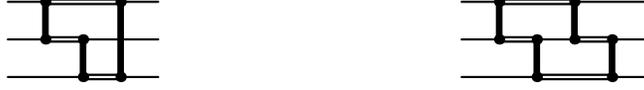

\centerline{
\includegraphics{tbd.1}
\hfil
\includegraphics{tbd.2}
}
\caption{If there are three consecutive pairs of pairwise distinct
  indices contained in $G|X$, there is always an $M$-alternating cycle
  for any $M\in \calM(G,X)$ using just the vertical edges (left); the
  same is true for any four consecutive pairs (right).  }
\label{fig:tbd}
\end{figure}

Consider the sequence $\Psi$ of 8 types occuring in $X$. 

Suppose $\Psi$ contains $1111$ as a subsequence, say
e.g. $\psi_2=\psi_3=\psi_4=\psi_5=1$. Then for $i=1$, $j=4$, $k=7$ the
condition (1) clearly holds. The condition (2) holds for $G|X_0$ in
this case, see Figure \ref{fig:psi}, left.
One can prove that the following subsequences 
are feasible as well by drawing triples of figures (we omit the details):
00011, 01011, 100000, 100010, 100100, 101000, 101010,
100110, 110110, 1000010, 1001010, 1010010, 00000000.

\begin{figure}
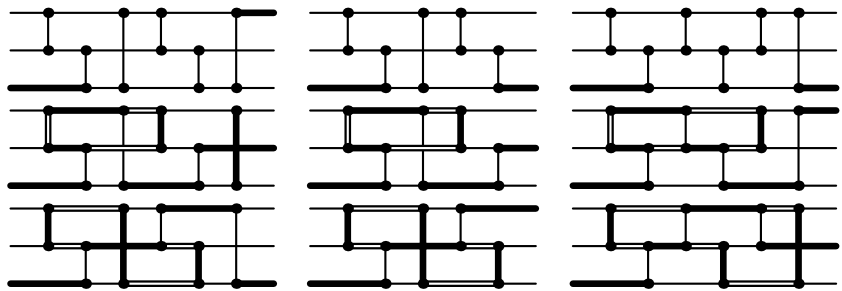

\centerline{
\begin{tabular}{ccc}
\includegraphics{psi2.3}
&
\includegraphics{psi1.1}
&
\includegraphics{psi3.1}
\\
\includegraphics{psi2.2}
&
\includegraphics{psi1.2}
&
\includegraphics{psi3.2}
\\
\includegraphics{psi2.1}
&
\includegraphics{psi1.3}
&
\includegraphics{psi3.3}
\end{tabular}
}
\caption{If $\Psi$ contains $1111$ (left), $111$ (center) or $0001$
  (right) as a subsequence, then for each perfect matching containing
  the bottom-leftmost edge such that $Y=2$ there is an alternating
  cycle. Observe that there is always one case out of three which is
  not possible.  }
\label{fig:psi}
\end{figure}

Similarly, if $\psi_1=\psi_2=\psi_3=1$, then we may pick $i=0$, $j=5$,
and $k=9$ (or even $k=8$ if $\psi_6=1$). In this case the condition
(2) holds for $G|X_1$, see Figure \ref{fig:psi}, center. Analogously
it works with 0001. It means $111**\,1$ and
$111****$ are feasible as well, and so are
$0001**\,1$ and $0001****$.

It remains to prove that $\Psi$ always contains at least one feasible
subsequence, which is a routine case analysis.
\end{proof}

\begin{lem}\label{l:2Burl} Let $(T,\phi)$ be a small-cut-decomposition
  of a cubic bridgeless graph $G$. Let $t_1,t_2 \in V(T)$ be a pair of
  adjacent vertices of degree $2$. Suppose that $|\phi^{-1}(f)|=2$ for
  every edge $f \in E(T)$ incident to $t_1$ or $t_2$.  Then
  $\phi^{-1}(\{t_1,t_2\})$ is a burl.
\end{lem}

\begin{proof}
  Let $t_0t_1t_2t_3$ be a subpath of $T$ and let
  $C_i=\phi^{-1}(t_{i-1}t_{i})$ for $i=1,2,3$ be an edge-cut of size
  $2$.  Assume that both $G|\phi^{-1}(t_1)$ and $G|\phi^{-1}(t_2)$
  have at most one perfect matching. By Lemma~\ref{l:4Burl} it
  suffices to show that $G|\phi^{-1}(\{t_1,t_2\})$ has at least two
  perfect matchings.  As the hub $G_1$ over $t_1$ is cubic and
  bridgeless it contains at least $2$ perfect matching avoiding any
  edge.  Let $e_1, e_2 \in E(G_1)$ be the edges in $E(G_1)-E(G)$
  corresponding to $C_1$- and $C_2$-contraction, respectively.  By
  assumption, at most one perfect matching of $G_1$ avoids both $e_1$
  and $e_2$.  It follows that either two perfect matchings of $G_1$
  avoid $e_1$ and contain $e_2$, or one avoids $e_1$ and $e_2$ and one
  avoids $e_1$ and contains $e_2$. Let $G_2$ be the hub over
  $t_2$. The symmetric statement holds for $G_2$. In any case, the
  perfect matchings in $G_1$ and $G_2$ can be combined to obtain at
  least two perfect matchings of $G|\phi^{-1}(\{t_1,t_2\})$.
\end{proof}

From the definition of a small-cut-decomposition, we immediately get
the following corollary:

\begin{cor}\label{c:2Burl}
  Let $(T,\phi)$ be a small-cut-decomposition of a cubic bridgeless
  graph $G$, and let $P$ be a path in $T$ in which every vertex has
  degree $2$.  Suppose there exist three edges $f_1$, $f_2$, $f_3$ of
  $T$ incident to vertices of $P$ such that
  $|\phi^{-1}(f_1)|=|\phi^{-1}(f_2)|=|\phi^{-1}(f_3)|=2$.  Then
  $\phi^{-1}(P)$ is a burl.
\end{cor}

Let $B_3$ denote the cubic graph consisting of two vertices joined by
three parallel edges.  Lemmas~\ref{l:3Burl} and~\ref{l:2Burl} imply
the following.

\begin{cor}\label{c:TreeBranch}
  Let $(T,\phi)$ be a small-cut-decomposition of a cubic bridgeless
  graph $G$ and let $P$ be a path in $T$ with $|V(P)|=32$. If for
  every $t \in V(P)$, $\deg_T(t)=2$ and the hub of $G$ at $t$ is
  isomorphic to $K_4$ or $B_3$, then $\phi^{-1}(P)$ is a burl.
\end{cor}

\begin{proof}
  If at least three edges incident to vertices in $V(P)$ correspond to
  edge-cuts of size $2$ in $G$ then the corollary holds by
  Corollary~\ref{c:2Burl}.  Otherwise, since there are $33$ edges of
  $T$ incident to vertices of $P$, there must be 11 consecutive edges
  incident to vertices in $P$ corresponding to edge-cuts of size $3$.
  In this case, the result follows from Lemma~\ref{l:3Burl}.
\end{proof}

\vskip 10pt
\begin{proof}[Proof of Lemma~\ref{l:Klee}]
  We proceed by induction on $|Z|$. If $|Z| \leq 6$ then $Z$ is a
  twig. In this case the lemma holds since $\beta_1 \geq \beta_2 +
  6\alpha$ by {\bf (\ref{e:const3})}. We assume for the remainder of
  the proof that $|Z| \geq 7$. It follows that $G'$ is not cyclically
  $4$-edge-connected, as $G'$ has no core. Therefore
  Corollary~\ref{c:ElemTwigs} is applicable to $G'$. Let $\calY$ be a
  collection of disjoint elementary twigs in $G'$ such that every
  relevant triangle in $G'$ is contained in an element of $\calY$, and
  let $(T,\phi)$ be a $\calY$-maximum small-cut decomposition of
  $G'$. By Corollary~\ref{c:ElemTwigs}, the hub at every $t \in V(T)$
  with $|\phi^{-1}(t)|\neq \emptyset$ is either an elementary twig, in
  which case $t$ is a leaf of $T$, or is cyclically
  $4$-edge-connected, in which case it is isomorphic to either $K_4$
  or $B_3$.
  
  In calculations below we will make use of the following claim: If
  $\deg_T(t)=2$ for some $t \in V(T)$, then $|\phi^{-1}(t)| \leq 2$.
  If this is not the case, the hub at $t$ is isomorphic to $K_4$, and
  at least three of its vertices must be vertices of $G$.  It follows
  that there is an edge $f\in E(T)$ incident to $t$ for which
  $|\phi^{-1}(f)|=2$. Let $v \in \phi^{-1}(t)$ be a vertex incident to
  an edge in $\phi^{-1}(f)$. Then $C=\phi^{-1}(f) \triangle \delta(v)$
  is a 3-edge-cut in $G$. As in the proof of
  Lemma~\ref{l:CutDecompose} we can split $t$ into two vertices
  $t',t''$ with $\phi^{-1}(t')=\{v\}$ and $\phi^{-1}(t'')=
  \phi^{-1}(t)\setminus v$. We now have $\phi^{-1}(t't'')=C$ and the
  new small-cut-decomposition contradicts the maximality of
  $(T,\phi)$. This completes the proof of the claim.

  Let $t_0 \in V(T)$ be such that $\phi^{-1}(t_0)$ contains the new
  vertex or one of the ends of the new edge in $G'$. Since $G$ is
  pruned, $G'$ contains at most one irrelevant triangle, and if such a
  triangle exists, at least one of its vertices lies in
  $\phi^{-1}(t_0)$. As a consequence, for any leaf $t\ne t_0$ of $T$,
  $\phi^{-1}(t)$ is a twig. Let $t^{*} \in V(T) \setminus \{t_0\}$ be
  such that $\deg_{T}(t^{*})\geq 3$ and, subject to this condition,
  the component of $T \setminus \{t^{*}\}$ containing $t_0$ is
  maximal. If $\deg_{T}(t) \leq 2$ for every $t \in V(T) \setminus
  \{t_0\}$, we take $t^{*}=t_0$ instead.

  Let $T_1,\ldots,T_k$ be all the components of $T \setminus
  \{t^{*}\}$ not containing $t_0$. By the choice of $t^{*}$, each
  $T_i$ is a path. If $|V(T_i)| \geq 33$ for some $1 \leq i \leq k$
  then let $T'$ be the subtree of $T_i$ containing a leaf of $T$ and
  exactly $32$ other vertices. Let $f$ be the unique edge in
  $\delta(T')$. 
  Let $H$ (resp.~$H'$) be the $\phi^{-1}(f)$-contraction of $G$ (resp.~$G'$)
containing $V(G')\setminus \phi^{-1}(T')$, and let $Z'$ consist of
$V(H') \cap Z$ together with the new vertex created by $\phi^{-1}(f)$-contraction (if it exists). If
$H$ is not pruned then it contains a unique irrelevant triangle
and we contract it, obtaining a pruned graph.
%
By the induction
  hypothesis, either $|Z'| \leq 6$ or we can find a foliage $\calX'$
  in $Z'$ with $\fw(\calX')\geq \alpha(|Z'|-2) + \beta_2$. If $|Z'|
  \leq 6$ let $\calX' := \emptyset$.

  Let $t'$ be a vertex of $T'$ which is not a leaf in $T$. Since
  $\deg_T(t')=2$, $|\phi^{-1}(t')|\neq \emptyset$. Therefore
  $\phi^{-1}(t')$ is isomorphic to $B_3$ or $K_4$ and we can apply
  Corollary~\ref{c:TreeBranch}. This implies that $\phi^{-1}(T')$
  contains an elementary twig and a burl that are vertex-disjoint,
  where the elementary twig is the preimage of the leaf. Further, we
  have $|\phi^{-1}(T')|\leq 8+ 2\cdot32 = 72$, since an elementary twig
  has size at most $8$ and the preimage of every non-leaf vertex of
  $T'$ has size at most $2$ by the claim above.  By
  Lemma~\ref{l:Contraction}(3), we can obtain a foliage $\calX$ in $Z$
  by adding the twig and the burl to $\calX'$ and possibly removing a
  burl (which can be a twig) containing the new element of $H'$
  created by $\phi^{-1}(f)$-contraction. It follows that if $|Z'| \geq
  7$ then
 $$\fw(\calX)\geq \alpha(|Z'|-2) + 2\beta_2 \geq (\alpha|Z|+\beta_2) -
 74\alpha +\beta_2 \geq \alpha|Z|+\beta_2,$$ by {\bf
   (\ref{e:const4})}, as desired. If $|Z'| \leq 6$ then $|Z| \leq 78$
 and
 $$\fw(\calX)\geq \beta_1 + \beta_2 \geq 78\alpha +\beta_2 \ge
 \alpha|Z| +\beta_2,$$ by {\bf (\ref{e:const5})}.

 It remains to consider the case when $|V(T_i)| \leq 32$ for every $1
 \leq i \leq k$. Suppose first that $t^{*} \neq t_0$ and that
 $|\phi^{-1}(T_0)| \geq 7$, where $T_0$ denotes the component of $T
 \setminus t^{*}$ containing $t_0$. Let $f_0$ be the edge incident to
 $t^*$ and a vertex of $T_0$. We form the graphs $H$, $H'$ and a set $Z'$
 by a $\phi^{-1}(f_0)$-contraction as in the previous case, and
 possibly contract a single irrelevant triangle. As before, we find a
 foliage $\calX'$ in $Z'$ with $\fw(\calX')\geq \alpha(|Z'|-2) +
 \beta_2$. Note that $\phi^{-1}(T_i)$ contains a twig for every $1\leq
 i\leq k$. By Lemma~\ref{l:Contraction}(3), we now obtain a foliage
 $\calX$ in $Z$ from $\calX'$, adding $k \geq 2$ twigs and possibly
 removing one burl (which can be a twig) from $\calX'$. We have
 $|\phi^{-1}(T_i)|\leq 8+31 \cdot 2=70$ for every $1 \leq i \leq k$,
 and $|\phi^{-1}(t^{*})| \leq 4$. Therefore $|Z| \leq |Z'|+70k+4$. It
 follows from {\bf (\ref{e:const5})} that
$$\fw(\calX) \geq \alpha(|Z'|-2) + \beta_2 + (k-1)\beta_1 \geq \alpha|Z|
+\beta_2 -76\alpha  +(k-1)(\beta_1-70\alpha) \geq \alpha|Z|+\beta_2.$$

Now we can assume $t^*=t_0$ or $|\phi^{-1}(T_0)|\leq 6$.  First
suppose $t^* \neq t_0$ but $|\phi^{-1}(T_0)|\leq 6$.  Then again
$|\phi^{-1}(t^*)|\leq 4$, so we have $|Z|\leq 70k+10$.  Let $\calX$ be
the foliage consisting of twigs in $T_1,\ldots,T_k$.  Thus by {\bf
  (\ref{e:const6})}, we have
\begin{equation*}\fw(\calX)=k\beta_1 \geq (\alpha|Z|+\beta_2) +
k(\beta_1-70\alpha)-10\alpha -\beta_2 \geq
\alpha|Z|+\beta_2.\end{equation*}

Finally we can assume $t^{*} = t_0$. Then $|\phi^{-1}(t^*)|\leq 4$,
unless $k=1$ and $\phi^{-1}(t^*)$ is an elementary twig.  In either
case, $|Z|\leq 70k+8$ and the equation above applies.
\end{proof}

\section{Proof of Lemma~\ref{l:MainInduction}}\label{s:MainInduction}

The following lemma is a direct consequence of a theorem of Kotzig,
stating that any graph with a unique perfect matching contains a bridge
(see~\cite{esperetkss10}).

\begin{lem}
\label{l:prop-double}
Every edge of a cyclically $4$-edge-connected cubic graph with at
least six vertices is contained in at least two perfect matchings.
\end{lem}

Let $G$ be a cubic graph. For a path $v_1v_2v_3v_4$, the graph
obtained from $G$ by {\em splitting along the path $v_1v_2v_3v_4$} is
the cubic graph $G'$ obtained as follows: remove the vertices $v_2$
and $v_3$ and add the edges $v_1v_4$ and $v'_1v'_4$ where $v'_1$ is
the neighbor of $v_2$ different from $v_1$ and $v_3$ and $v'_4$ is the
neighbor of $v_3$ different from $v_2$ and $v_4$. The idea of this
construction (and its application to the problem of counting perfect
matchings) originally appeared in~\cite{voorhoeve79}. We say that a
perfect matching $M$ of $G$ is a \emph{canonical extension} of a
perfect matching $M'$ of $G'$ if $M \triangle M' \subseteq E(G)
\triangle E(G')$, i.e. $M$ and $M'$ agree on the edges shared by $G$
and $G'$.

\begin{lem}
\label{l:splitting}
Let $G$ be a cyclically 4-edge-connected cubic graph with $|V(G)| \geq
6$. If $G'$ is the graph obtained from $G$ by splitting along some
path $v_1v_2v_3v_4$, then
\begin{enumerate}
\item $G'$ is cubic and bridgeless;
\item $G'$ contains at most $2$ irrelevant triangles;
  \item $\fw(G) \geq \fw(G')-2\beta_1$;
  \item Every perfect matching $M'$ of $G'$ avoiding the edge $v_1v_4$
    has a canonical extension in $G$.
\end{enumerate}
\end{lem}
\begin{proof}\mbox{}\smallskip
\begin{enumerate}

\item The statement is a consequence of an easy lemma
  in~\cite{esperetkk11}, stating that the cyclic edge-connectivity can
  drop by at most two after a splitting.

\item Since $G$ is cyclically 4-edge-connected and has at least six
  vertices, it does not contain any triangle. The only way an
  irrelevant triangle can appear in $G'$ is that $v_1$ and $v_4$ (or
  $v_1'$ and $v_4'$) have precisely one common neighbor (if they have
  two common neighbors, the two arising triangles share the new edge
  $v_1v_4$ or $v_1'v_4'$ and hence, are relevant).

\item At most two burls from a foliage of $G'$ contain $\{v_1,v_4\}$
  or $\{v_1',v_4'\}$. Therefore, a foliage of $G$ can be obtained from
  any foliage of $G'$ by removing at most two burls (observe that this
  is precisely here that we use the fact that being a burl is a local
  property, independent of the rest of the graph).

\item The canonical extension is obtained (uniquely) from $M' \cap
  E(G)$ by adding either $v_2v_3$ if $v'_1v'_4 \not\in M'$ or
  $\{v'_1v_2,v_3v_4'\}$ if $v'_1v'_4 \in M'$.\qedhere
\end{enumerate}
\end{proof}

\vskip 10pt
\begin{proof}[Proof of Lemma~\ref{l:MainInduction}]
  We proceed by induction on $|V(G)|$. The base case $|V(G)|=6$
  holds by Lemma~\ref{l:prop-double} and {\bf (\ref{e:const7})}.

  For the induction step, consider first the case that $G$ is
  cyclically 4-edge-connected. Fix an edge $e=uv \in E(G)$. Our goal
  is to show that $e$ is contained in at least
  $2^{\alpha|V(G)|-\fw(G)+\gamma}$ perfect matchings.

  Let $w \neq u$ be a neighbor of $v$ and let $w_1$ and $w_2$ be the
  two other neighbors of $w$. Let $x_i,y_i$ be the neighbors of $w_i$
  distinct from $w$ for $i=1,2$. Let $G_1,\ldots,G_4$ be the graphs
  obtained from $G$ by splitting along the paths $vww_1x_1$,
  $vww_1y_1$, $vww_2x_2$ and $vww_2y_2$. Let $G'_i$ be obtained from
  $G_i$ by contracting irrelevant triangles for $i=1,\ldots,4$. By
  Lemma~\ref{l:splitting}(2) we have $|V(G'_i)| \geq |V(G)|-6$.

  Suppose first that one of the resulting graphs, without loss of
  generality $G'_1$, does not have a core. By Corollary~\ref{c:Klee},
  Lemma~\ref{l:Triangle} and Lemma~\ref{l:splitting}, we have $$\alpha
  |V(G)| \leq \alpha (|V(G'_1)| + 6) \leq \fw(G'_1)+ 7\alpha -\beta_2
  \leq \fw(G_1)+ 7\alpha -\beta_2 \leq \fw(G)+ 2\beta_1 +7\alpha
  -\beta_2.$$ Therefore $$\alpha|V(G)|-\fw(G)+\gamma \leq \gamma +
  2\beta_1 +7\alpha -\beta_2 \leq 1$$ by {\bf (\ref{e:const8})} and the
  lemma follows from Lemma~\ref{l:prop-double}.

  We now assume that all four graphs $G'_1,\ldots,G'_4$ have a
  core. By Lemma~\ref{l:splitting}(4), every perfect matching
  containing $e$ in $G_i$ canonically extends to a perfect matching
  containing $e$ in $G$. Let $S$ be the sum of the number of perfect
  matchings of $G_i$ containing $e$, for $i\in\{1,2,3,4\}$. By
  induction hypothesis and Lemmas~\ref{l:Triangle}
  and~\ref{l:splitting}, $S \ge 4 \cdot
  2^{\alpha(|V(G)|-6)-\fw(G)-2\beta_1+\gamma}$. On the other hand, a
  perfect matching $M$ of $G$ containing $e$ is the canonical
  extension of a perfect matching containing $e$ in precisely three of
  the graphs $G_i$, $i\in\{1,2,3,4\}$. For instance if $w_1y_1,ww_2
  \in M$, then $G_2$ is the only graph (among the four) that does not
  have a perfect matching $M'$ that canonically extends to $M$ (see
  Figure~\ref{f:voo}). As a consequence, there are precisely $S/3$
  perfect matchings containing $e$ in $G$. Therefore,
$$m^{\star}(G) \ge \tfrac{4}{3}\cdot
2^{\alpha(|V(G)|-6)-\fw(G)-2\beta_1+\gamma} \geq
2^{\alpha|V(G)|-\fw(G)+\gamma},$$ by {\bf (\ref{e:const9})}, as
desired.

\begin{figure}[htbp]
\centering
\includegraphics[scale=0.35]{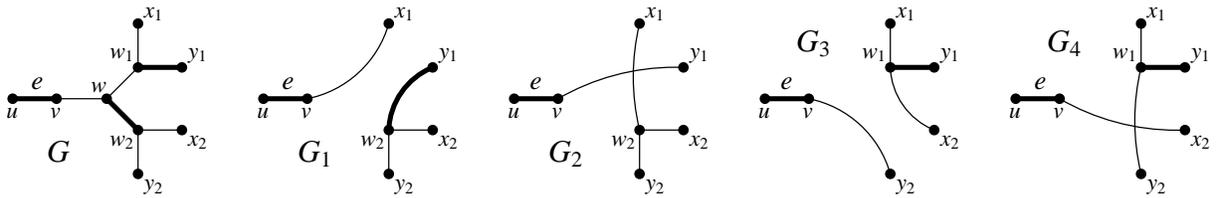}
\caption{Perfect matchings in only three of the $G_i$'s canonically
  extend to a given perfect matching of $G$ containing $e$.} \label{f:voo}
\end{figure}

  It remains to consider the case when $G$ contains a cyclic edge-cut
  $C$ of size at most $3$. Suppose first that for such edge-cut $C$,
  both $C$-contractions $H_1$ and $H_2$ have a core. Then, by
  Lemma~\ref{l:Contraction}(3), $\fw(G) \geq \fw(H_1)+\fw(H_2)-2\beta_1$
  and, by induction hypothesis, applied to $H_1$ and $H_2$ (after
  possibly contracting one irrelevant triangle in each) and
  Lemma~\ref{l:Contraction},
$$m^{\star}(G) \geq m^{\star}(H_1)m^{\star}(H_2) \geq
2^{\alpha|V(G)|-4\alpha-\fw(G)-2\beta_1 +2\gamma} \geq
2^{\alpha|V(G)|-\fw(G)+\gamma},$$ by {\bf (\ref{e:const10})}, as desired.
Finally, if for every cyclic edge-cut $C$ of size at most $3$ only one
$C$-contraction has a core, we apply Corollary~\ref{c:ElemTwigs} to
$G$. Let $(T,\phi)$ be the resulting small-cut-decomposition of $G$.
There exists a unique vertex $t \in V(T)$ such that the hub $H$ of $G$
at $t$ is cyclically $4$-edge-connected with $|V(H)|\geq 6$. Let
$T_1,\ldots, T_k$ be the components of $T - t$ and let $Z_i =
\phi^{-1}(V(T_i))$. We apply Lemma~\ref{l:Klee} to
$Z_1,\ldots,Z_k$. Note that Lemma~\ref{l:Klee} is indeed applicable,
as $G$ is pruned, and therefore every triangle in $G$ belongs to an
elementary twig. Consequently, no edge-cut corresponding to an edge of $(T,
\phi)$ separates exactly $3$ vertices of $G$.

Let $\calX_1,\calX_2,\ldots,\calX_k$ be the foliages satisfying the
lemma. Let $\calX_0$ be the maximal foliage in $H$ avoiding new
vertices and edges created by contraction of the edge-cuts
$\delta(Z_1),\ldots,\delta(Z_k)$. Then $\fw(\calX_0) \geq
\fw(H)-k\beta_2$, as $H$ contains no twigs (it is cyclically
4-edge-connected). Since $\calX_0 \cup \calX_1 \cup \ldots \cup \calX_k$
is a foliage in $G$ we have
$$\fw(G) \geq \fw(H)-k\beta_2 + \sum_{i=1}^{k}\fw(\calX_i) \geq \fw(H) +
\alpha \sum_{i=1}^{k}|Z_i|,$$ by the choice of
$\calX_1,\ldots,\calX_k$.  It remains to observe that
$$m^{\star}(G) \geq m^{\star}(H) \geq 2^{\alpha|V(H)|-\fw(H)+\gamma}
\geq 2^{\alpha(|V(G)|-\sum_{i=1}^{k}|Z_i|)-\fw(H)+\gamma} \geq
2^{\alpha|V(G)|-\fw(G)+\gamma},$$ by the above.
\end{proof}

\section{Concluding remarks}
\subsection{Improving the bound.} We expect that the bound in
Theorem~\ref{t:Main} can be improved at the expense of more careful
case analysis.  In particular, it is possible to improve the bound on
the length of the path in Corollary~\ref{c:TreeBranch}. We have chosen
not to do so in an attempt to keep the argument as straightforward as
possible. 

In~\cite{cyganps09} it is shown that for some constant $c>0$ and every
integer $n $ there exists a cubic bridgeless graph on at least $n$
vertices with at most $c2^{n/17.285}$ perfect matchings.

\subsection{Number of perfect matchings in $k$-regular graphs.}
In~\cite[Conjecture 8.1.8]{lovaszp86} the following generalization of
the conjecture considered in this paper is stated. A graph is said to
be \emph{matching-covered} if every edge of it belongs to a perfect
matching.

\begin{conj}\label{conj:kregular}
  For $k \geq 3$ there exist constants $c_1(k), c_2(k) > 0$ such that
  every $k$-regular matching-covered graph contains at least $c_2(k)
  c_1(k)^{|V(G)|}$ perfect matchings. Furthermore, $c_1(k) \to \infty$
  as $k \to \infty$.
\end{conj}

While our proof does not seem to extend to the proof of this
conjecture, the following weaker statement can be deduced from
Theorem~\ref{t:Main}. We are grateful to Paul Seymour for suggesting
the idea of the following proof.

\begin{thm}\label{t:kregular}
  Let $G$ be a $k$-regular $(k-1)$-edge-connected graph on $n$
  vertices for some $k \geq 4$. Then
$$\log_2{m(G)} \geq (1-\tfrac1k)(1-\tfrac2k) \tfrac{n}{\ceps}.$$
\end{thm}
\begin{proof}
  It follows by Edmonds' characterization of the perfect matching
  polytope~\cite{edmonds65} that there exists a probability
  distribution $\bM$ on $\calM(G)$ such that for every edge $e\in
  E(G)$, $\Pr[e\in \bM]=\tfrac1k$.

  We choose a triple of perfect matchings of $G$ as follows. Let $M_1
  \in \calM(G)$ be arbitrary.  We have $$\Ex[|\bM \cap
  M_1|]=\frac{n}{2k}.$$ Therefore we can choose $M_2 \in \calM(G)$ so
  that $|M_2 \cap M_1|\leq \frac{n}{2k}$. Let $Z \subseteq V(G)$ be
  the set of vertices not incident with an edge of $M_1 \cap
  M_2$. Then $|Z| \geq (1-\tfrac1k)\, n$. For each $v \in Z$ we
  have $$\Pr[\bM \cap \delta(\{v\}) \cap (M_1 \cup M_2) =
  \emptyset]=1-\tfrac2k.$$ Therefore the expected number of vertices
  whose three incident edges are in $\bM$, $M_1$ and $M_2$
  respectively, is at least $(1-\tfrac1k)(1-\tfrac2k)\,n$. It follows
  that we can choose $M_3 \in \calM(G)$ so that the subgraph $G'$ of
  $G$ with $E(G')=M_1 \cup M_2 \cup M_3$ has at least
  $(1-\tfrac1k)(1-\tfrac2k)\,n$ vertices of degree three. Note that
  $G'$ is by definition matching-covered. It follows that the only
  bridges in $G'$ are edges joining pairs of vertices of degree
  one. Let $G''$ be obtained from $G'$ by deleting vertices of degree
  one and replacing by an edge every maximal path in which all the
  internal vertices have degree two. The graph $G''$ is cubic and
  bridgeless and therefore by Theorem~\ref{t:Main} we have
$$\log_2{m(G)}> \log_2{m(G')} \geq \log_2{m(G'')} \geq
\tfrac1{\ceps}|V(G'')| \geq (1-\tfrac1k)(1-\tfrac2k)
\tfrac{n}{\ceps},$$ as desired.\end{proof}

\bibliographystyle{plain}

\end{document}